\documentclass[12pt,letterpaper]{article}
\usepackage{t1enc}
\usepackage[latin1]{inputenc}
\usepackage[english]{babel}
\begin{document}

\title{On the period function of Newtonian systems}

\author{A. Raouf Chouikha \footnote{Universite Paris 13 LAGA, Villetaneuse 93430, chouikha@math.univ-paris13.fr},  Mohsen Timoumi \footnote{Faculté des Sciences de Monastir, Tunisie, m\_timoumi@yahoo.com}}

\date{}
\maketitle {}

\begin{abstract}
  We study the existence of centers of planar autonomous system of the form
  $$(S) \quad \dot x=y,\qquad \dot y = -h(x) - g(x)y - f(x)y^2.$$
  We are interested in the period function $T$ around a center $0$. A sufficient condition for the isochronicity of $(S)$ at $0$ is given. Such a condition is also necessary when $f,g,h$ are analytic functions. In that case a characterization of isochronous centers of  system $(S)$ is given. Some applications will be derived. In particular, new families of isochronous centers will be described.\footnote
{2000 Mathematics Subject Classification  \ 34C15, 34C23, 34C25, 34C37, 37G15 }. 
\end{abstract}

\section {Introduction}
The study of the generalized Lienard equations of the form
\begin{equation}\ddot x+f(x)\dot x^2+g(x)\dot x+h(x)=0\end{equation}
or, its equivalent two-dimensional form,
\begin{equation} \quad \dot x=y,\qquad \dot y = -h(x) - g(x)y - f(x)y^2,\end{equation}
(where $\dot \xi = \frac{d\xi}{dt}$ and $\ddot \xi = \frac{d^2\xi}{dt^2}$)
holds an important place in the theory of dynamical systems. This equation is sometime called in the litterature as Langmiur equation [8]. 
The Langmuir equation governs the space-charge current in an electron tube is a special case, as is the equation for the brachistochrone, [8]. Others models can be reduced to system (2) by means of suitable transformations.\\

In this paper, we are interested in the conditions under which (2) has a center as well as a center of
constant period (or, alternatively, for which (1) has a non-isolated periodic solution
with locally constant period). Such centers are called isochronous and it was shown
by Poincare that this implies that the critical
point is locally linearizable (without time scaling) and vice versa. 
Interest in the phenomena of isochronicity arises from several areas. One aspect,
is the desire to understand what mechanisms allow systems to be locally linearizable
in the polynomial case. Another aspect is to understand better the nature of
the period function for a family of closed orbits. This in turn is useful in the study
of critical periods and bifurcations for
planar systems.\\ 
We know very few significant classes of system classified to date: Quadratic systems were classified by Loud 
as well as some cubic and quartic systems. We refer to [2] for a survey concerning these problems.\\ The main tool of this paper is to derive necessary and sufficient conditions in order to system (2) has a center at the origin well as to produce  necessary and sufficient conditions for this center to be isochronous. These results extend the ones of [1], [4] and [5]. \\

Notice that Volokitin and Ivanov [11] have already interested in existence of isochronous centers for system (S) and proved that in the polynomial case system (S) have not a non trivial isochronous center. More precisely, they proved when $f,g$ and $h$ are polynomials then (S) cannot commute with any polynomial system nonproportional to it. Therefore we suppose in the sequel that $f,g$ and $h$ are analytic non polynomial functions.\\  

On the other hand, Cherkas ([2]) introduced the following center condition for the classical Lienard system

\begin{equation} \dot x=y,\qquad \dot y = -h(x) - g(x)y,\end{equation}

\bigskip

{\bf Proposition 1}\quad  {\it  Assume that $g$ and $h$ are analytic with $g(0)=h(0)=0$ and $h'(0)=1$. Then the lienard system (3) has a nondegenerate center at the origin if and only if $G(x) = \int_0^x g(\xi) d\xi$ can be expressed as an analytic function of $H(x) = \int_0^x h(\xi) d\xi$}

\bigskip

In other words, this means $G(x) = \phi (H(x))$ where $\phi$ is some analytic function with $\phi (0)=0$. The preceding result has been improved by Christopher [6] who obtained a global necessary and sufficient condition when $g(x)$ and $h(x)$ are polynomials of some polynomial $p(x)$ of degree greater than one. But in general it is not easy to find $p(x)$.\\ 

Concerning the classification of isochronous centers of Lienard systems Christopher and Devlin [7] proved the following characterization

\bigskip

{\bf Proposition 2}\quad  {\it   System $$ \dot x=y,\qquad \dot y = -h(x) - g(x)y,$$ with $h$ and $g$ analytic
functions of $x$ such that $h(0) = g(0) = 0, h'(0) = 1$ and $x h(x)>0$ in a neighborhood of the origin has an isochronous center at the origin if and only if 
$$ h(x)=s(x)s'(x)\biggl(1+{1\over
s^4(x)}\biggl(\int_0^xs(\xi)g(\xi)\,d\xi\biggr)^2\biggr),$$ 
where the analytic function $s(x)$ solves the functional equation} 
$$G(x - 2s(x)) = G(x),  \qquad G(x) = \int_0^x g(\xi) d\xi, \quad s(0) = 0, s'(0) = 1. $$ 

\bigskip

In particular,
 when in addition $h(x)$ is odd, then
(3) has an isochronous center at the origin  if and only if
$g(x)$ is odd and
$$
h(x) = x + {1\over x^3}\biggl(\int_0^x \xi g(\xi)\,d\xi\biggr)^2.
$$

\bigskip

This latter result which means $ s(x) \equiv x$ was first proved by Sabatini [9] for
Li\'enard systems that are not necessarily analytic.

The proof of this theorem follows from bringing the system (3)
into a normal form using function derived from the complex
separatrices of the system at the origin.
In [3] an alternative easier proof of Proposition 1 is presented. \\

More recently, Amel'kin ([1], Theorems 2 and 4) gave some improvments of the last result and proved the following general characterization

\bigskip

{\bf Proposition 3}\quad  {\it  System $$\dot x=y,\qquad \dot y = -h(x) - g(x)y,$$ with $h$ and $g$ analytic
functions of $x$ such that $h(0) = g(0) = 0, h'(0) = 1$ has an isochronous center at the origin if and only if $g(x)$ is an odd function and the relation $$h(x) = x + \frac{1}{x^3}[\int_0^x \xi g(\xi) d\xi]^2$$ holds.}

\bigskip

{\bf Proposition 4}\quad  {\it Let us write $h(x)=x+ \sum_{i\geq 2} a_i x^i, \ g(x) = \sum_{j\geq 2} b_j x^i.$ Then, system 
$$\dot x=y,\qquad \dot y = -h(x) - g(x)y,$$ 
has an isochronous center at the origin $0$ if and only if} 
$$a_{2k}= 0, \quad a_{2k+1}= \sum_{0\leq m\leq k-1} (\frac{b_{2m+1}}{2m+3})(\frac{b_{2k-2m-1}}{2k-2m-1}), \quad k = 1,2,...$$

\bigskip

Let us return to system (S) and define the integral $F(x) = \int_0^x f(\xi) d\xi.$
 For the special case $ g(x)\equiv 0$\ in (2) one proved [4] the following which concerns the Lienard type system\\

{\bf Proposition 5}\quad {\it Let $f,h$ be analytic function in a neighborhood $N_0$ of $0$  
and $x h(x) > 0$ for $x \neq 0$ \ 
then System  \begin{equation} \quad \dot x=y,\qquad \dot y = -h(x)  - f(x)y^2,\end{equation} has an isochronous center at the origin $0$ if and only if 
$$\frac {X}{1+u(X)} = h(x) e^{F(x)}$$ where \ $X$ \ is defined by\ $\frac {1}{2} X^2 = \int_0^x h(s) e^{2F(s)} ds$ \ and \ $u(X)$ \ is an odd function such that \ $ \phi(x) = \int_0^x e^{F(s)} ds = X + \int_0^X u(t) dt$\ and \ $\frac {X}{\phi (x)} > 0$.\\
In particular, when \ $u(X) \equiv 0$ \  then \ $0$\ is an isochronous center if and only if \ $h(x) = e^{-F(x)}\phi (x)$.}\\

  \bigskip
  
 Concerning the center condition of system (2) we derive the following 
 
 \bigskip

{\bf Theorem A}\quad {\it Assume that $f, g$ and $h$ are analytic with $f(0)=g(0)=h(0)=0$ and $h'(0)=1$. Then the system (2) has a nondegenerate center at the origin if and only if $ \int_0^x g(\xi)e^{F(\xi)} d\xi $ can be expressed as an analytic function of $\int_0^x h(\xi) e^{2F(\xi)} d\xi.$}\\
 \bigskip
 
In this paper we prove the following results

\bigskip

{\bf Theorem B}\quad {\it  Let $f,g$ and $h$ analytic
functions such that $h(0) = f(0) = 0, h'(0) = 1$ and $x h(x)>0$ in a neighborhood of the origin. Then, system (2)  has an isochronous center at the origin if and only if 
$$(C)\qquad  h(x) e^{F(x)}=s(y)s'(y)\biggl(1+{1\over
s^4(y)}\biggl(\int_0^ys(\xi)g(x(\xi))\,d\xi\biggr)^2\biggr),$$ 
where $F(x) = \int_0^x f(t) dt, \ y = \int_0^x e^{F(\xi)}d\xi$ and $s(y)$ solves the functional equation 
$$\tilde G(y - 2s(y)) = \tilde G(y),  \quad \tilde G(y) = \int_0^y g(u(\xi)) d\xi = \int_0^x g(\xi)) e^{F(\xi))}d\xi, \quad s(0) = 0, s'(0) = 1. $$  
\\
In particular, when $g \equiv 0$ then condition (C) reduces to} $$s(y) = \sqrt{2\int_0^y h(x(\xi)) e^{F(x(\xi))}d\xi} = \sqrt{2\int_0^x h(\xi) e^{2F(\xi)}d\xi}. $$\\

\bigskip

{\bf Theorem C}\quad {\it Let $f,g$ and $h$ analytic
functions of $x$ such that $f(x)$ is an odd function and $h(0) = f(0) = 0, h'(0) = 1$. Then, system (2)  has an isochronous center at the origin if and only if $g$ is an odd function and the relation $$ h(x) e^{F(x)}= y + {1\over y^3}\biggl(\int_0^y \xi g(x(\xi))\,d\xi\biggr)^2
$$ holds 
where \ $y = \int_0^x e^{F(\xi)}d\xi$  and $F(x)= \int_0^x f(\xi) d\xi$.}\\

\section{Proofs}

 Let us consider the diffeomorphic change $x = u(y) $ in equation (1)
$$\ddot x+f(x)\dot x^2+g(x)\dot x+h(x)=0$$ such that $u(0)=0, u'(0)=1$ with
$$\frac {dy}{dx}= e^{F(x)}, \quad where \quad F(x) = \int_0^x f(\xi) d\xi.$$ So, $y=u^{-1}(x)=\int_0^x e^{F(\xi)}d\xi.$\\ 
The derivative of $u$ must verify $$\frac {du}{dy}= e^{-F(u(y))}.$$ 

{\bf Lemma 3}\quad The function $u$ defined by $\frac {du}{dy}= e^{-F(u(y))}, u(0)=0, u'(0)=1$ is an analytic diffeomorphism and equation (1) is equivalent to 
\begin{equation}
  \ddot y+g(u(y))\dot y + h(u(y))e^{F(u(y))}=0.
\end{equation} 
\bigskip

Indeed, since $ \frac {du}{dy}(0)= 1$ then $u$ is well defined and it is an analytic diffeomorphism. Moreover, we get $\dot x = u'(y)\dot y$ and 
$\ddot x = u''(y) \dot y^2 + u'(y) \ddot y$. \\ Replacing in equation (1) we then obtain 
$$u'(y) \ddot y + [f(u(y))u'^2(y)+u''(y) ]\dot y^2 + g(u(y)u'(y)\dot y + h(u(y)) = 0.$$
From $\frac {dy}{dx}= e^{F(x)}$ we deduce $f(u(y))u'^2(y)+u''(y) =0$ and then equation (1) is equivalent to
\begin{equation}
  \ddot y+g(u(y))\dot y + \frac {h(u(y))}{u'(y)}=0.
\end{equation}

Finally equation (6) is a Lienard equation $$\ddot y + \tilde g(y) \dot y + \tilde h(y) = 0$$
where $\tilde g(y) = g(u(y))$ and $\tilde h(y)= h(u(y))e^{F(u(y))}.$ On the others words, $\tilde g(y) = g(x)$ and $\tilde h(y)= h(x)e^{F(x)}.$\\
 Moreover, by Proposition 1 equation (6) has a nondegenerate center at the origin if and only if $\tilde G(y) = \int_0^y \tilde g(\xi)d\xi$ can be expressed as an analytic function of $\tilde H(y) = \int_0^y h(\xi) d\xi$. Or equivalently $ \int_0^x g(\xi)e^{F(\xi)} d\xi $ can be expressed as an analytic function of $\int_0^x h(\xi) e^{2F(\xi)} d\xi$.\\

 By Proposition 2 equation (6) admits an isochronous center at $0$ if and only if \begin{equation} \tilde h(y)=s(y)s'(y)\biggl(1+{1\over
s^4(y)}\biggl(\int_0^y s(\xi)\tilde g(\xi)\,d\xi\biggr)^2\biggr),\end{equation}  
where $s(y)\neq y$ solves the functional equation 
$$\tilde G(y - 2s(y)) = \tilde G(y),  \qquad \tilde G(y) = \int_0^y \tilde g(\xi) d\xi, \qquad s(0)=0, s'(0)=1.$$
Moreover, this may also be written $$ \tilde G(y) = \tilde G(u^{-1}(x)) = \int_0^x  g(\xi) e^{F(\xi)}d\xi.$$

That means by Proposition 2 
$$ h(u(y))e^{F(u(y))} = s(y)s'(y) \biggl(1+{1\over
s^4(y)}\biggl(\int_0^ys(\xi) g(u(\xi))\,d\xi\biggr)^2\biggr).$$ Or equivalently
$$h(x)e^{F(x)}= s(u^{-1}(x))s'(u^{-1}(x)) \biggl(1+{1\over
s^4(u^{-1}(x))}\biggl(\int_0^xs(u^{-1}(x)) g(\xi)e^{F(\xi)}\,d\xi\biggr)^2\biggr)$$ where $y = u^{-1}(x)=\int_0^x e^{F(\xi)}d\xi.$ \\ Thus, Theorem B is proved.

Moreover, applying Proposition 3 we then obtain that $$(5) \qquad \ddot y+g(u(y))\dot y + h(u(y))e^{F(u(y))} =0$$  has an isochronous center at the origin if and only if $\tilde g(y)=g(u(y))$ is an odd function and the relation $$\tilde h(y)= h(u(y))e^{F(u(y))} = y + \frac{1}{y^3}[\int_0^y \xi g(\xi) d\xi]$$ holds where \ $y = \int_0^x e^{F(\xi)}d\xi$. Then this  achieves the proof of Theorem C.\\

\section{Consequences}

{\bf Corollary 1}\quad {\it Let $f,g$ and $h$ analytic
functions of $x$ such that $h(0) = f(0) = 0, h'(0) = 1$ and $x h(x)>0$ in a neighborhood of the origin. Then, system (2)  has an isochronous center at the origin if and only if 
$$(C')\qquad h(x)e^{2F(x)}=\sigma(x)\sigma'(x)\biggl(1+{1\over
\sigma^4(x)}\biggl(\int_0^x\sigma(\xi)g(\xi)e^{F(x)}\,d\xi\biggr)^2\biggr),$$ 
where $\sigma(x)\neq x$ solves the functional equation }
$$\tilde G(u^{-1}(x) - 2\sigma(x)) = \tilde G(u^{-1}(x))= \int_0^x g(\xi))e^{F(\xi)} d\xi, \quad \sigma(0) = 0, \sigma'(0) = 1. $$ 

\bigskip

Indeed let us consider the function $$\sigma = s\circ u^{-1}, \quad  i.e. \quad \sigma (x) = s(u^{-1}(x)) = s(y).$$ Then, $$\sigma'(x) = \frac{d\sigma (x)}{dx}=\frac{ds}{dy}\frac{du^{-1}}{dx} = \frac{ds}{dy} e^{F(x)}$$ and $\frac{ds}{dy}(y) = \sigma '(x) e^{-F(x)}.$\\ Thus, the condition (C) of Theorem B \\
$ h(x) e^{F(x)}=s(y)s'(y)\biggl(1+{1\over
s^4(y)}\biggl(\int_0^ys(\xi)g(x(\xi))\,d\xi\biggr)^2\biggr)$ \ 
may be written  
$$(C')\qquad h(x)e^{2F(x)}=\sigma(x)\sigma'(x)\biggl(1+{1\over
\sigma^4(x)}\biggl(\int_0^x\sigma(\xi)g(\xi)e^{F(x)}\,d\xi\biggr)^2\biggr),$$ 
where $\sigma(x)\neq x$ solves the functional equation 
$$\tilde G(u^{-1}(x) - 2\sigma(x)) = \tilde G(u^{-1}(x))= \int_0^x g(\xi))e^{F(\xi)} d\xi, \quad \sigma(0) = 0, \sigma'(0) = 1. $$ 
 
Turn now to the case when $g\equiv 0$. Equality (7) becomes \\ $h(u(y)) e^{F(u(y))} = s(y) s'(y)$ which easily implies \\ $\int_0^y h(u(\xi)) e^{F(u(\xi))}d\xi = \frac {1}{2} s^2(y)$ since $s(0)=0.$ In fact, this means by Proposition 3 : \ $s(y)= X.$

\bigskip

As corollaries we then get Proposition 2 when $ f(x)\equiv 0$ and Proposition 4 when $ g(x)\equiv 0$. Thus, Theorems B and C generalize previous results.\\

\bigskip

{\bf Corollary 2}\quad {\it Let $f,g$ and $h$ analytic
functions of $x$ such that $f(x)$ and $g(x)$ or $h(x)$ is an odd function such that $h(0) = f(0) = 0, h'(0) = 1$ and $x h(x)>0$ in a neighborhood of the origin. Then, system (2) has an isochronous center at the origin if and only if 
$$ h(x) e^{F(x)}= y + {1\over y^3}\biggl(\int_0^y \xi g(u(\xi))\,d\xi\biggr)^2
$$
where} \ $y = \int_0^x e^{F(\xi)}d\xi$  and $F(x)= \int_0^x f(\xi) d\xi$.

\bigskip 

Indeed, when $g$ is odd means obviously $s(y) = y$ since $\tilde G(-y) = \tilde G(y)$ and $G(y - 2s(y)) = \tilde G(y)$. Or equivalently by Corollary 1 $\sigma (x) = y = u^{-1}(x).$ Moreover, it has been proven (see [6] for example) that when the Lienard system has an isochronous center at $0$ and $g$ is odd then $h$ is also necessarily odd.

\bigskip 

{\bf Corollary 3}\quad {\it  Let $f,g,h$ analytic functions of $x$ and $g$ or $h$ is an odd function such that $h(0) = f(0) = 0, h'(0) = 1$ and $x h(x)>0$ in a neighborhood of the origin. Then, system (2) has an isochronous center at the origin if and only if} $$h(x)f(x) + \frac{dh}{dx}(x) = 1 + \frac{2 g(u(y))}{y^2} \int_0^y \xi g(u(\xi))\,d\xi - \frac{1}{y^4}\biggl(\int_0^y \xi g(u(\xi))\,d\xi\biggr)^2  .$$

\bigskip 
Moreover, as a direct consequence of Theorem C and Corollary 3, we find again the result proved in [11].

\bigskip 
{\bf Corollary 4}\quad {\it Let $f,g,h$ polynomials verifying $h(0) = f(0) = 0, h'(0) = 1$ and $x h(x)>0$ in a neighborhood of the origin. Then system (2) has no isochronous center at the origin $0$.}

\bigskip 

\section {Some applications} 
In this section thanks to {\it Maple} we present examples in order to illustrate our results. First, it is interesting to note that one can find  examples of systems (2) with different functions $f, g, h$ satisfying Theorem C, thus ensuring the isochroniciy of the origin $0$. We have of course to looking for examples of non polynomial functions. What is amazing what one can also find other more general examples where $f$ and $h$ are depending of an arbitrary function $g$.

\subsection {A rational case} Consider at first the following simple case 
$$f(x)= 2 g(x) = \frac{2x}{x^2+1}$$ and 
$$h(x)= \frac{x+\frac{x^3}{3}+\frac{9x^3}{(3+x^2)^2}(1+\frac{x^2}{5})^2}{x^2+1}.$$
Then, $e^{F(x)}=x^2+1$ and $y=x+\frac{x^3}{3}$ of equation (6). It implies $$\tilde h(y) = h(x) e^{F(x)} = y + \frac{1}{y^3}[\int_0^y \xi g(\xi) d\xi]^2= x+\frac{x^3}{3}+\frac{9x^3}{(3+x^2)^2}(1+\frac{x^2}{5})^2.$$ 
 One deduces that equation $ \ddot y + \tilde g(y) \dot y + \tilde h(y) = 0$\ has an isochronous center at $0$ where $$\tilde g(y)= g(u(y))=g(x),\quad \tilde h(y)= h(u(y)) e^{F(u(y))}=h(x) e^{F(x)}.$$  
We thus prove that the system $$\dot x = -y ,\qquad \dot y = (\frac{2x}{x^2+1})y^2 -(\frac{x}{x^2+1})y +  \frac{x+\frac{x^3}{3}+\frac{9x^3}{(3+x^2)^2}(1+\frac{x^2}{5})^2}{x^2+1}$$ has an isochronous center at the origin $0$. Notice that one has necessary $\mid x\mid < 1$ and the derivative $g'(x)= \frac {1-x^2}{(1+x^2)^2}$ is positive.\\ Thus, one gets an example of isochronous center when $f,g,h$ are rational functions. We hope there exist many others and to classify all these centers for the rational case.

\subsection {A family of isochronous centers} Now look at somewhat more general cases, where hypothesis $g'(x) > 0$ still to be necessary. 
Consider again equation (6) $$ \ddot y + \tilde g(y) \dot y + \tilde h(y) = 0$$ where $\tilde g(y)= g(u(y))=g(x),\ \tilde h(y)= h(u(y)) e^{F(u(y))}=h(x) e^{F(x)}$ and $y=\int_0^x e^{F(s)}ds.$ Suppose $\tilde g(y) = y$ or equivalently $g'(x)= e^{F(x)}$. By Proposition 3 equation (6) has an isochronous center if and ony if $\tilde g(y)$ is odd and 
$$\tilde h(y)= h(u(y)e^{F(u(y))} = y + \frac{1}{y^3}[\int_0^y \xi g(\xi) d\xi].$$ Then  
$$\tilde h(y) = y + \frac{1}{y^3}[\int_0^y \xi^2  d\xi]= y + \frac{1}{9}y^3.$$ It implies $$h(x) e^{F(x)}= h(x) g'(x) = g(x) + \frac{1}{9}g^3(x).$$ Thus $h(x)= \frac{g(x)}{g'(x)}+ \frac{g^3(x)}{9g'(x)} $ and $ f(x)= \frac{g''(x)}{g'(x)}$ where $g'(x)$ is positive.\\
Finally, replacing in equation (6) we then proved the following

\bigskip  

{\bf Proposition 6}\quad {\it System $$\dot x = -y ,\qquad \dot y = \frac{g''(x)}{g'(x)} y^2 - g(x) y + \frac{g(x)}{g'(x)}+ \frac{g^3(x)}{9g'(x)}$$ has an isochronous center at the origin if and only if $g(x)$ is an odd analytic function such that $g'(x) > 0$.}

\bigskip 
This proposition generalizes the special case $g(x) = x $ which yields the wellknown Lienard polynomial system $$\dot x = -y ,\qquad \dot y = xy + \frac{9 x+x^3}{9}.$$ 

\subsection {Others examples} By the same way we may find many other systems of this type by different choice of coefficients of system (2).\\
Let us consider for equation (1) the following case $$y= \sinh (g(x)).$$
Then \ $y' = \cosh (g(x)) g'(x) =  e^{F(x)}$\ implying $g'(x) > 0.$ It follows $$\int _0^y s \tilde g(s)ds = \int_0^x \sinh (g(x)) g(x) e^{F(x)} dx = \frac{1}{4}g(x) \cosh(2g(x)) - \frac{1}{8} \sinh(2g(x)).$$ It implies $$h(x)e^{F(x)} = \sinh(g(x)) + \frac{1}{(\sinh(g(x))^3} [\frac{1}{4}g(x) \cosh(2g(x)) - \frac{1}{8} \sinh(2g(x))]^2$$ and $f(x)= \tanh(g(x)) g'(x).$ We thus prove

\bigskip  

{\bf Proposition 7}\quad {\it System $$\dot x = -y ,$$ $$ \dot y = \tanh(g(x)) g'(x) y^2 - g(x) y + \sinh(g(x)) + \frac{1}{(\sinh(g(x))^3} [\frac{1}{4}g(x) \cosh(2g(x)) - \frac{1}{8} \sinh(2g(x))]^2$$ has an isochronous center at the origin if and only if $g(x)$ is an odd analytic function such that $g'(x) > 0$.} 

\bigskip   

By the same manner we may prove

\bigskip  

{\bf Proposition 8}\quad {\it System $$\dot x = -y ,$$ $$ \dot y = [\frac{g''(x)}{g'(x)}-\frac{g(x) g'(x)}{1+g^2(x)}] y^2 - g(x) y + \sinh^{-1}(g(x)) + \frac{1}{(\sinh^{-1}(g(x))^3} [\sinh^{-1}(g(x))\sqrt{1+g^2(x)}-g(x)]^2$$ has an isochronous center at the origin if and only if $g(x)$ is an odd analytic function such that $g'(x) > 0$.} 

\bigskip

\section {An analytic involution}

In this section we present results of smaller importance but related to what precedes. We know that concerning periodic orbits the invariance of turning points by the energy brings up an involution. It is proposed to clarify the relationship between this involution and isochronous centers for potential and Lienard systems

\subsection {For the potential equation}
 
Let us consider the potential equation \begin{equation} \ddot x + h(x) = 0\end{equation} where $h$ is analytic. Let $$H(x)= \int h(t) dt$$ be the potential. We suppose in the sequel $h(0)=0, h'(0) > 0$ and $xh(x) > 0$ so that the origin $0$ is a center of (8). 
 Since the potential $H(x)$ has a local minimum at $0$, then we may consider an involution $A$ by 
$$ H(A(x)) = H(x) \ and \ A(x) x < 0 $$ for all $x \in [a,b]$ and $x \neq 0$. So, any closed orbit is $A$-invariant and $A$ exchanges the turning points: $b=A(a)$. \\ In fact, $A(x)$ is well defined in the interval $[a,b]$. To see that, set the function $$\rho(x)=\frac{x-A(x)}{2}.$$ This function is such that $\rho(A(x)) = - \rho(x)$ and $\rho'(x) = \frac{1-A'(x)}{2}$. Since $A'(x)<0$ we get $\rho'(x) >0$ and therefore $\rho$ is an analytic diffeomorphism on $[a,b]$. Then $A(x) = \rho^{-1}(-\rho(x))$ is well defined.  \\
The following holds

 \bigskip

{\bf Proposition 9}\quad {\it Suppose $0$ is a center of $(8)$, then the analytic potential $H$ is $A$-invariant and under the action of the involution $A$ and equation $(8)$ is equivalent to the following lienard type equation 
$$ \frac{d^2 A(y)}{dy^2} \dot y^2 + \frac{d A(y)}{dy} \ddot y + h(A(y)) = 0. \qquad (8') $$ Moreover, the following are equivalent : \\ 
1) - $0$ is an isochronous center of $(8)$ \\ 
2) - $0$ is an isochronous center of $(8')$ \\
3) - $h(x)= \frac{1}{4} (x - A(x))(1 - \frac{d A(x)}{dx})$ \\
4) - $A(x) = x - 2 \sqrt{2 H(x)}$}

 \bigskip
 
 {\bf Proof}\quad Indeed the involution $A$ is such that $x = A(y)$ is a diffeomorphism verifying $A(0)=0,  \frac{d A(y)}{dy}(0)=-1$.\\
 The period function of equation  (6) is $$ T(c) = \sqrt{2}\int_a^b \frac{dx}{\sqrt{c-H(x)}}$$ where $c$ is the energy level such that $H(a) = H(b) =c.$ Since the involution $A$ exchanges the turning points ( $A(a) = b$) the period may also be written 
 $$ T(c) = - \sqrt{2}\int_a^b \frac{A'(y) dy}{\sqrt{c-H(y)}}.$$
 Thus the involution $A$ must verify the following identity
 $$\int_a^b \frac{[1+A'(y)] dy}{\sqrt{c-H(y)}}=0.$$

Let us write $\tilde f(y) = \frac{A''(y)}{A'(y)}$ and $\tilde h(y) = \frac {h(A(y))}{A'(y)}$ then $(8')$ may becomes
$$ \ddot y + \tilde f(y) \dot y^2 + \tilde h(y) = 0. \qquad (8'')$$ By proposition 3 equation $(8'')$ has an isochronous center if and only if $$\tilde h(y) e^{\tilde F(y)}= \frac {X}{1+u(X)}$$ where $\frac {1}{2} X^2 = \int_0^x \tilde h(s) e^{2\tilde F(s)} ds$ \ and \ $u(X)$ \ is an odd function such that \ $ \phi(x) = \int_0^x e^{\tilde F(s)} ds = X + \int_0^X u(t) dt.$\ We then get obviously $$ \tilde h(y) e^{\tilde F(y)}= \frac {h(A(y))}{A'(y)} A'(y) = h(A(y)).$$ This means conditions 1) and 2) are equivalent. \\
Conditions $1)$ and $4)$ are equivalent by Proposition 3-1 of [5].\\ Deriving $A(x) = x - 2 \sqrt{2 H(x)}$ one gets $A'(x) = 1 - \frac{2 h(x)}{2 \sqrt{2 H(x)}}$ which implies $$(1-A'(x))2 \sqrt{2 H(x)} = 2 h(x).$$ Replacing $2 \sqrt{2 H(x)}$ by $ \frac{x-A(x)}{2}$ one then obtains\\ condition  $3)$ $h(x)= \frac{1}{4} (x - A(x))(1 - \frac{d A(x)}{dx})$ and conversely.\\

\subsection {For the Lienard equation} 

Consider again the Lienard equation \begin{equation} \ddot x+g(x)\dot x+h(x)=0\end{equation} or equivalently the Lienard system
 $$ \dot x=y,\qquad \dot y = -h(x) - g(x)y. $$ We suppose $g(0)=h(0)=0$ \ and\ $h'(0)>0.$ Then the integral
 $$H(x) = \int_0^x g(\xi) d\xi $$ has a minimum local at $0$. We then define an involution $A$ by
 $$H(A(x))= H(x) \quad and \quad A(x) x < 0.$$
 We prove the following
 
 \bigskip

{\bf Proposition 10}\quad {\it Suppose $0$ is a center of $(9)$, then the analytic potential $H$ is $A$-invariant and under the action of the involution $A$ and equation $(9)$ is equivalent to the following Langmiur type equation 
$$ \frac{d^2 A(y)}{dy^2} \dot y^2 + \frac{d A(y)}{dy} \ddot y + g(A(y))\frac{d A(y)}{dy} \dot y + h(A(y)) = 0.  \qquad (9')$$ Moreover, the following are equivalent :} 
\begin{center} 
{\it 1) - $0$ is an isochronous center of $(9)$ \\ 
2) - $0$ is an isochronous center of $(9')$ 
 $${\it 3) -} \quad h(x)=s(x)s'(x)\biggl(1+{1\over
s^4}\biggl(\int_0^xs(\xi)g(\xi)\,d\xi\biggr)^2\biggr),$$ 
where $s(x)=\frac{1}{2}[x- A(x)]$ 
 $${\it 4) -} \quad h(A(y))=\sigma(y)\sigma'(y)\biggl(1+{1\over
\sigma^4}\biggl(\int_0^y\sigma(\xi)g(A(\xi))A'(\xi)\,d\xi\biggr)^2\biggr),$$ 
where $\sigma(y)=\frac{1}{2}[y- A(y)].$}
\end{center}

 \bigskip
 {\bf Proof} \quad By Proposition 2 we get $1)  \Leftrightarrow 3) $. By Corollary 1 we get $2)  \Leftrightarrow 4) $ since 
 $e^{F(y)} = A'(y)$. Moreover, from $\sigma = s\circ A,$ or $ s = \sigma \circ A$ we then deduce $3)  \Leftrightarrow 4) -$.

 \newpage

REFERENCES

\bigskip

[1] V.V. Amel'kin, {\it Strong isochronicity of the Lienard system}, Diff. Equations, vol. 42, n 5, 615-618, (2006).

\smallskip

[2] J. Chavarriga, M. Sabatini ,{\it A survey of isochronous centers},
Qual. Theory of Dyn. Systems vol 1 , 1-70, (1999).

\smallskip

[3]\ A. R. Chouikha, {\it Monotonicity of the period function for some planar differential
systems. Part II : Lienard and related systems},
Applicationes Mathematicae, 32 no. 4, 405-424, (2005).

\smallskip

[4]\ A. R. Chouikha, {\it Isochronous centers of Lienard type equations and applications},
J. Math. Anal. Appl. 331, 358-376 (2007).

\smallskip

[5]\ A. R. Chouikha,  {\it Period function and characterization of isochronous potentials}\quad arXiv:1109.4611, (2011).

\smallskip
[6] \ C. Christopher, {\it An algebraic approach to the classification of centers in polynomial Lienard systems}, J. Math. Anal. Appl., 229, n 1, 319-329 (1999).

\smallskip
[7] \ C. Christopher and J. Devlin, 
             {\it On the classification of Lienard Systems with Amplitude-Independent Periods},  J. Diff. Eq., {\bf 200}, 1-17, (2004).
\smallskip

[8]\ H. Denman and L. Buch, {\it Variational treatment for equation of the Langmuir type }, SIAM, J. Appl. Math., vol 40, 2, 279-282, (1981).           

\smallskip
[9] \ M. Sabatini,  {\it On the period function of Lienard systems,} \quad J. of Diff. Eq., 152,  467-487, (1999).

\smallskip

[10]\ M. Sabatini, {\it  On the period function of} $\ddot x + f(x)\dot x^2 + g(x) = 0$ ,\\
J. Diff. Eq. 196, no. 1,  151-168 (2004).

\smallskip
[11] E.P. Volokitin and V.V. Ivanov,  {\it Isochronicity and commutation of polynomial vector fields},\
Siberian Math. J., 40,1, 23-38 (1999).

\end{document}